\documentclass[a4paper,reqno,12pt]{amsart}
\usepackage[latin2]{inputenc}
\usepackage[T1]{fontenc}
\usepackage{amsmath}
\usepackage{amsthm}
\usepackage{amsfonts}
\usepackage{amssymb}

\theoremstyle{plain}
\newtheorem*{KKthm}{Kaloujnine--Krasner Theorem}
\newtheorem{thm}{Theorem}[section]

\newtheorem{prop}[thm]{Proposition}

\theoremstyle{definition}

\usepackage[pdftex]{graphicx}

\DeclareMathOperator{\Ker}{Ker}
\newenvironment{biz}{\rm \trivlist\item[\hskip \labelsep{\it
Proof.\quad}]}{\hfill\qed\par\medskip\endtrivlist}

\def\act#1#2{{^{#1}\kern -2pt {#2}}}
\def\acth#1#2{{^{#1}\kern -1pt {#2}}}

\begin{document}

\title{On extensions of completely simple semigroups by groups}
\author{Tam\'as D\'ek\'any}
\address{Bolyai Institute, University of Szeged, Aradi v\'ertan\'uk tere 1, Szeged, Hungary, H-6720; fax: +36 62 544548}
\email{tdekany@math.u-szeged.hu}

\date{August 10, 2013}

\thanks{Research supported by the Hungarian National Foundation for Scientific
Research grant no.\ K083219, K104251,
and by the European Union, cofunded by the European Social Fund, under the project 
no.\ TÁMOP-4.2.2.A-11/1/KONV-2012-0073.
}

\begin{abstract}
An example of an extension of a completely simple semigroup $U$ by a group $H$ is given which cannot be  
embedded into the wreath product of $U$ by $H$. 
On the other hand, every central extension of $U$ by $H$ is shown to be embeddable in the wreath product of $U$ by $H$, and 
any extension of $U$ by $H$ is proved to be embeddable in a semidirect product of 
a completely simple semigroup $V$ by $H$ where the maximal subgroups of $V$ are direct powers of those of $U$.
\end{abstract}

\maketitle

\section{Introduction}

Group extensions play a fundamental role both in the structure theory and in the theory of varieties of groups. 
The following basic result was proved in 1950, see \cite{KK}:

\begin{KKthm}
Any extension of a group $N$ by a group $H$ is embeddable in the wreath product of $N$ by $H$.
\end{KKthm}

Note that the wreath product of $N$ by $H$ is a special semidirect product of a direct power of $N$ by $H$, see the details below.

Completely simple semigroups are often considered as a natural and close generalization of groups. 
In this paper we establish that the Kaloujnine--Krasner Theorem
fails to hold for extensions of completely simple semigroups by groups, 
where extension is taken within the class of regular, or equivalently, of completely simple semigroups. 
However, it is valid within the class of central completely simple semigroups, and
a slightly weaker embedding theorem holds in general: each extension of a completely simple semigroup $U$ 
by a group $H$ is embeddable in a semidirect product of 
a completely simple semigroup $V$ by $H$ where $V$ is ``close'' to $U$ 
(e.g., the maximal subgroups of $V$ are direct powers of those of $U$), and in the special case where $U$ is a group,
the embedding given in the proof coincides with that well known from a standard proof of the Kaloujnine--Krasner Theorem.

\section{Preliminaries}

For the undefined notions and notation the reader is referred to \cite{Ho} and \cite{PeRe}.

Let $T$ be a semigroup and let $H$ be a group with identity element $N$. 
We say that $H$ {\em acts on $T$ (on the left by automorphisms)}
if a map $H \times T\to T$, $(A,t) \mapsto \acth {A}{t}$ is given such that, 
for all $t,u \in T$ and $A,B \in H$, we have
$$ \act Nt = t,\ \act{BA}t = \act B{(\acth At)}\hbox{ and }\acth A{(tu)}=\acth At\cdot\act Au.$$ 
If $T$ is a semigroup acted upon by the group $H$ then an associative  multiplication can be defined on the set 
$T\times H$ by the rule
$(t,A)(u,B)=(t\cdot \act{A}{u},A\cdot B)$.
The semigroup so obtained is called the {\em semidirect product of $T$ by $H$ (with respect to the given action)},
and is denoted by $T\rtimes H$.
Observe that if $T$ is also a group then this construction is the usual semidirect product of groups.

Notice that the second projection $T\rtimes H$ induces a group congruence $\rho$ on $T\rtimes H$, and 
$(T\rtimes H)/\rho$ is isomorphic to $H$. 
Moreover, the identity element of the group $(T\rtimes H)/\rho$ is
$$\Ker \rho=\{(t,N)\in T\rtimes H : t\in T\},$$
called the {\em kernel of $\rho$}, which is a subsemigroup in $T\rtimes H$ isomorphic to $T$.

Note that a semidirect product of a completely simple semigroup by a group is completely simple.
Green's relations $\mathcal{L}$ and $\mathcal{R}$ on a semidirect product can be described in the following way.

\begin{prop}\label{RL}
Let $T\rtimes H$ be a semidirect product of a regular semigroup $T$ by a group $H$.
Then, for any elements $(t,A),(u,B)$, we have
\begin{itemize}
\item[(1)] $(t,A)\,\mathcal{R}\,(u,B)$ in $T\rtimes H$ if and only if $t\,\mathcal{R}\, u$ in $T$,
\item[(2)] $(t,A)\,\mathcal{L}\,(u,B)$ in $T\rtimes H$ if and only if $\act{A^{-1}}{t}\,\mathcal{L}\, \act{B^{-1}}{u}$ in $T$.
\end{itemize}
\end{prop}

Now let $T$ be a semigroup and let $H$ be a group. 
An action of $H$ on the direct power $T^H$ can be defined in the following natural way:
for any $f\in T^H$ and $A\in H$, let $\act Af$ be the element of $T^H$ where $B(\act{A}{f})=(BA)f$
for any $B\in H$. 
The semidirect product $T^H\rtimes H$ defined in this way is called the {\em wreath product of $T$ by $H$}, and is denoted by $T\wr H$. 

By the Rees--Suschkewitsch Theorem, a completely simple semigroup is in most cases represented throughout the paper as a 
Rees matrix semigroup with a normalized sandwich matrix.
A completely simple semigroup is called {\em central} if the product of any two of its idempotents lies in the centre 
of the containing maximal subgroup. 
It is well known that a Rees matrix semigroup $\mathcal{M}[G;I,\Lambda;P]$ with $P$ normalized is central if and only if
each entry of $P$ belongs to the centre of $G$.

The group congruences of a Rees matrix semigroup with a normalized sandwich matrix are characterized as follows.
 
\begin{prop}\label{kongruencia}
Let $S=\mathcal{M}[G;I,\Lambda;P]$ be a Rees matrix semigroup where $P$ is normalized. 
Assume that $N$ is a normal subgroup of $G$ such that every entry of $P$ belongs to $N$. 
Define a relation $\rho$ on $S$ such that, for every $(i,g,\lambda)$, $(j,h,\mu)\in S$, let
\begin{equation*}
(i,g,\lambda)\,\rho\, (j,h,\mu) \quad \hbox{if and only if}\quad gh^{-1}\in N.
\end{equation*}
Then $\rho$ is a group congruence on $S$ such that $S/\rho$ is isomorphic to $G/N$ and $\Ker \rho=\mathcal{M}[N;I,\Lambda;P]$.

Conversely, every group congruence on $S$ is of this form for some normal subgroup $N$ of $G$ where
all entries of $P$ belong to $N$.
\end{prop}

This proposition implies that the kernel of any group congruence of a completely simple semigroup is completely simple. 
Conversely, it is routine to check that if $S$ is a regular semigroup and $\rho$ is a group congruence on $S$
such that $\Ker \rho$ is a completely simple subsemigroup of $S$
then $S$ is necessarily completely simple.
This allows us to extend the notion of a group extension in the following manner:
given completely simple semigroups $S,U$ and a group $H$, 
we say that {\em $S$ is an extension of $U$ by $H$} if there exists a group congruence $\rho$ on $S$ such that 
$S/\rho$ is isomorphic to $H$ and $\Ker \rho$ is isomorphic to $U$. 
For example, a semidirect product of a completely simple semigroup $T$ by a group $H$ is an extension of $T$ by $H$. 

Now we present an isomorphic copy of the wreath product $T\wr H$ of a Rees matrix semigroup $T=\mathcal{M}[G;I,\Lambda;P]$ by a group $H$
which allows us to make the calculation in the next section easier.
First, it is routine to see that the direct power $T^H$ is isomorphic to $\mathcal{M}[G^H;I^H,\Lambda^H;P^H]$
where $P^H=(p_{\xi \eta}^H)$ is the following sandwich matrix: 
for any $\xi\in \Lambda^H$ and $\eta\in I^H$ we have $Ap_{\xi \eta}^H=p_{A\xi,A\eta}\ (A\in H)$.
Moreover, the action in the definition of the wreath product determines the following action when replacing
$T^H$ by $\mathcal{M}[G^H;I^H,\Lambda^H;P^H]$: 
for any $A\in H$ and $(\eta,f,\xi)\in \mathcal{M}[G^H;I^H,\Lambda^H;P^H]$ we have 
$\acth{A}{(\eta,f,\xi)}=(\acth{A}{\eta},\act{A}{f},\acth{A}{\xi})$, 
where $\acth{A}{\eta}\in I^H$, $\act{A}{f}\in G^H$ and $\acth{A}{\xi}\in \Lambda^H$ 
are the maps defined by $B(\acth{A}{\eta})=(BA)\eta$, $B(\act{A}{f})=(BA)f$ and $B(\acth{A}{\xi})=(BA)\xi$, respectively, for every $B\in H$.

Notice that, for any $A\in H$, we have
\[ A(\act{B}{p_{\xi\eta}^H}) = (AB)p_{\xi\eta}^H = p_{(AB)\xi,(AB)\eta} = 
p_{A(\act{B}{\xi}),A(\act{B}{\eta})} = Ap_{\act{B}{\xi}, \act{B}{\eta}}^H, \]
and so
\begin{equation}\label{M}
\act{B}{p_{\xi\eta}^H}=p_{\act{B}{\xi}, \act{B}{\eta}}^H
\end{equation}
for any $B\in H$.

Finally, we sketch a standard proof of the Kaloujnine--Krasner Theorem. 

Let $G$ be an extension of $N$ by $H$. 
Without loss of generality, we can assume that $N$ is a normal subgroup of $G$ and $H=G/N$. 
Choose and fix an element $r_A$ from each coset $A$ of $N$ in $G$ such that $r_N$ is the identity element of $G$. 
It is straightforward to check that the map
\begin{equation}\label{KKemb}
\kappa\colon G\rightarrow N^H\rtimes H,\ g\mapsto (f_g,gN)\ \hbox{where}\ 
f_g\colon H\rightarrow N,\ A\mapsto r_Agr_{A\cdot gN}^{-1} 
\end{equation}
is an embedding. Since $\kappa$ is a morphism, the equality
\begin{equation}\label{hatas}
f_{gh}=f_g\cdot \act{gN}f_{h}
\end{equation}
holds for every $g,h\in G$.

\section{Embeddability in a wreath product}

In this section first we notice that the Kaloujnine--Krasner Theorem can be easily extended to central completely simple semigroups.
Moreover, we establish that it fails in general: we present a completely simple semigroup which is an extension of
a completely simple semigroup by a group, and is not embeddable in their wreath product.

Let $S=\mathcal{M}[G;I,\Lambda;P]$ be an extension of a completely simple semigroup $U$ by a group $H$ where $P$ is chosen to be normalized. 
By Proposition \ref{kongruencia}, we can assume that there is a normal subgroup $N$ of the group $G$ 
such that all entries of the sandwich matrix $P$ belong to $N$, and we have $H=G/N$ and $U=\mathcal{M}[N;I,\Lambda;P]\subseteq S$.

First suppose that $S$ is central, i.e.,  each entry of $P$ belongs to the centre of the group $G$.
Note that, in this case, $U$ is necessarily also central. 
In this case, we can mimic the proof of the Kaloujnine--Krasner Theorem sketched in the previous section. 
For, it is routine to check by applying \eqref{KKemb} and \eqref{hatas} that the map
\[ \nu\colon S\rightarrow U\wr H=U^H\rtimes H,\ (i,g,\lambda)\mapsto (f_g^{i\lambda},gN) \]
where
\[ f_g^{i\lambda}\colon H\rightarrow U,\ A\mapsto(i,Af_g,\lambda) \]
is an embedding. This verifies the following statement.

\begin{prop}
Each central completely simple semigroup which is an extension of a (necessarily also central) completely simple semigroup 
$U$ by a group $H$ is embeddable in the wreath product of $U$ by $H$.
\end{prop}

Now we turn to investigating the general case where $S$ is an arbitrary completely simple semigroup. 
Suppose that there exists an embedding $S\rightarrow U\wr H$, i.e.,\ an embedding
\begin{equation}\label{beagyazas}
\varphi\colon S\rightarrow \mathcal{M}[N^H;I^H,\Lambda^H;P^H]\rtimes H 
\end{equation}
where $\mathcal{M}[N^H;I^H,\Lambda^H;P^H]\rtimes H$ is the isomorphic copy of $U\wr H$ introduced in the previous section.
Proposition \ref{RL} implies that, in the semidirect product $\mathcal{M}[N^H;I^H,\Lambda^H;P^H]\rtimes H$, 
we have 
$$[(\eta_1,f_1,\xi_1),A]\,\mathcal{R} \,[(\eta_2,f_2,\xi_2),B]$$
if and only if 
$(\eta_1,f_1,\xi_1)\,\mathcal{R}\,(\eta_2,f_2,\xi_2)$ in $\mathcal{M}[N^H;I^H,\Lambda^H;P^H]$, 
and this is the case if and only if $\eta_1=\eta_2$. 
Moreover, 
$$[(\eta_1,f_1,\xi_1),A]\,\mathcal{L}\,[(\eta_2,f_2,\xi_2),B]$$
if and only if 
$\act{A^{-1}}{(\eta_1,f_1,\xi_1)}\,\mathcal{L}\,\act{B^{-1}}{(\eta_2,f_2,\xi_2)}$ in $\mathcal{M}[N^H;I^H,\Lambda^H;P^H]$, 
and this is the case if and only if $\act{A^{-1}}{\xi_1}=\act{B^{-1}}{\xi_2}$. 
Thus we see that the $\mathcal{R}$-class of an element $[(\eta,f,\xi),A]$ depends only on $\eta$, 
while its $\mathcal{L}$-class depends only on $\xi$ and $A$. 
Since the morphism $\varphi$ sends $\mathcal{R}$-equivalent elements to $\mathcal{R}$-equivalent elements, 
and $\mathcal{L}$-equivalent elements to $\mathcal{L}$-equivalent elements, we obtain that,
for each $i\in I$, there exists $\eta_i\in I^H$, and for each $(A,\lambda)\in H\times \Lambda$, there exists $\xi_{A,\lambda}\in \Lambda^H $, 
such that, for every $g\in G$, we have
\[ (i,g,\lambda)\varphi=[(\eta_i,f_g^{i\lambda},\xi_{gN,\lambda}),gN] \]
for some $f_g^{i \lambda}\in N^H$.

Since $\varphi$ is a morphism, we have
\[ [(\eta_i,f_g^{i\lambda},\xi_{gN,\lambda}),gN][(\eta_j,f_h^{j\mu},\xi_{hN,\mu}),hN] = 
[(\eta_i,f_{gp_{\lambda j}h}^{i\mu},\xi_{ghN,\mu}),ghN] \]
for any $i,j\in I$, $g,h\in G$ and $\lambda,\mu\in \Lambda$. 
This equality holds if and only if
\begin{equation}\label{A}
f_g^{i\lambda}\cdot p_{\xi_{gN,\lambda},\act{gN}{\eta_j}}^H\cdot \act{gN}{f_h^{j\mu}} = f_{gp_{\lambda j}h}^{i\mu}
\end{equation}
for any $i,j\in I$, $g,h\in G$ and $\lambda,\mu\in \Lambda$, and
\begin{equation}\label{B}
\act{gN}{\xi_{hN,\mu}} = \xi_{ghN,\mu}
\end{equation}
for any $g,h\in G$ and $\mu \in \Lambda$. 
Notice that \eqref{B} is equivalent to requiring that
\[ \xi_{A,\mu}=\acth{A}{\xi_{N,\mu}} \]
for every $\mu\in \Lambda$ and $A\in H$.
Therefore, later on, we shortly write $\xi_\mu$ and $\act{A}{\xi_\mu}$ instead of $\xi_{N,\mu}$ and $\xi_{A,\mu}$, respectively.

By \eqref{M}, equality \eqref{A} is equivalent to
\begin{equation}\label{C}
f_{gp_{\lambda j}h}^{i\mu} = f_g^{i\lambda}\cdot \act{gN}{p_{\xi_\lambda\eta_j}^H}\cdot \act{gN}{f_h^{j\mu}}.
\end{equation}
Substituting $g=h=1$ and  $g=p_{\lambda i}^{-1}$, $h=1$, $j=i$, respectively, 
where $1$ denotes the identity element of $N$, we obtain from \eqref{C} that
\begin{equation}\label{D}
f_{p_{\lambda j}}^{i\mu} = f_1^{i\lambda} p_{\xi_\lambda\eta_j}^H f_1^{j\mu},
\end{equation}
\begin{equation*}
f_1^{i\mu}=f_{p_{\lambda i}^{-1}p_{\lambda i}1}^{i\mu}=f_{p_{\lambda i}^{-1}}^{i\lambda} p_{\xi_\lambda\eta_i}^H f_1^{i\mu},
\end{equation*}
and the latter implies
\begin{equation}\label{E}
f_{p_{\lambda i}^{-1}}^{i\lambda} = (p_{\xi_\lambda\eta_i}^H)^{-1}.
\end{equation}

If $p_{\lambda i}=1$ then the map
\begin{equation}\label{iota}
\iota_{i \lambda}\colon G\rightarrow N^H\rtimes H,\quad g\mapsto (p_{\xi_\lambda\eta_i}^H f_g^{i \lambda},gN)
\end{equation}
is an injective group morphism. 
For, it is injective since $\varphi$ is injective, and by \eqref{C}, we have 
$p_{\xi_\lambda\eta_i}^H f_{gh}^{i \lambda} = p_{\xi_\lambda\eta_i}^H f_{gp_{\lambda i}h}^{i \lambda} = 
p_{\xi_\lambda\eta_i}^H f_g^{i \lambda}\cdot \act{gN}{(p_{\xi_\lambda\eta_i} f_{h}^{i \lambda})}$, 
and so
\[ (p_{\xi_\lambda\eta_i}^H f_g^{i \lambda},gN)(p_{\xi_\lambda\eta_i}^H f_h^{i \lambda},hN) = 
(p_{\xi_\lambda\eta_i}^H f_{gh}^{i \lambda},ghN).\]

We now give a suitable group $G$, a normal subgroup $N$ of $G$ and a Rees matrix semigroup $S=\mathcal{M}[G;I,\Lambda;P]$ 
for which no such injective morphism $\varphi$ exists.

Let $G$ be the non-commutative group of order $21$. 
To ease our calculations, we present $G$ in the form $G=\mathbb{Z}_7\rtimes \left[ \overline{2}\right]$ 
where $\mathbb{Z}_7$ is the additive group of the ring of residues modulo $7$, 
$\left[ \overline{2}\right]=\{\overline{1},\overline{2},\overline{4}\}$ is the subgroup of the (multiplicative) group 
of units of the same ring generated by $\overline{2}$, and $\left[ \overline{2}\right]$ acts on $\mathbb{Z}_7$ by multiplication. 
The second projection of $G$ is a morphism onto $[\overline{2}]$, 
its kernel is $N=\{(a,\overline{1}):a\in \mathbb{Z}_7\}$ isomorphic to $\mathbb{Z}_7$, 
and the map $H=G/N\rightarrow [\overline{2}]$, $(a,k)N\mapsto k$ is an isomorphism. 
For our later convenience, we identify $H$ with $[\overline{2}]$ via this isomorphism. 
Let $I=\Lambda=\{1,2\}$, and denote by $P$ the normalized sandwich matrix of type $\Lambda\times I$ over $G$ consisting of the elements 
$p_{11}=p_{12}=p_{21}=(\overline{0},\overline{1})$, the identity element of $N$, and $p_{22}=(\overline{1},\overline{1})\in N$, 
an element of order $7$.

Now we assume that $\varphi$ is an embedding of the form \eqref{beagyazas} from this Rees matrix semigroup $S$, 
and apply the general properties deduced so far for this $S$.

The elements of order $3$ in $G$ and $N^H\rtimes H$ play crucial role in our argument. 
Observe that $(\overline{0},\overline{2})$ and $(\overline{0},\overline{4})$ are mutual inverse elements of $G$ of order $3$. 
Moreover all the elements of order $3$ in $N^H\rtimes H$ are of the form $(t,\overline{2})$ or $(t,\overline{4})$. 
Let us mention, although we do not need it explicitly, that $(t,\overline{2})$ and $(t,\overline{4})$ are of order $3$ if and only if 
$(\overline{1}t)\cdot(\overline{2}t)\cdot(\overline{4}t)=(\overline{0},\overline{1})$.

Applying the injective group morphism $\iota_{11}\colon G\rightarrow N^H\rtimes H$ defined in \eqref{iota}, 
we see that $p_{22}\iota_{11}=(h,\overline{1})$ with $h=p_{\xi_1\eta_1}^H f_{p_{22}}^{11}$. 
Since the image of an element of order $3$ has order $3$, the following two cases occur:

{\it Case 1\/}: $(\overline{0},\overline{4})\iota_{11}=(t,\overline{4})$. 
Then we obtain $((\overline{0},\overline{4})^{-1}p_{22}(\overline{0},\overline{4}))\iota_{11}=$ $
(\act{\overline{2}}{(t^{-1})},\overline{2})(h,\overline{1})(t,\overline{4})=
(\act{\overline{2}}{(t^{-1})} \cdot \act{\overline{2}}{h} \cdot \act{\overline{2}}{t},\overline{1})=
(\act{\overline{2}}{h},\overline{1})$. 
On the other hand, $(\overline{0},\overline{4})^{-1}p_{22}(\overline{0},\overline{4}) =
(\overline{0},\overline{2})(\overline{1},\overline{1})(\overline{0},\overline{4}) =(\overline{2},\overline{1})=
(\overline{1},\overline{1})^2=p_{22}^2$, 
and so $p_{22}^2\iota_{11}= (h,\overline{1})(h,\overline{1})=(h^2,\overline{1})$. 
Thus $\act{\overline{2}}{h}=h^2$ which implies, for any $a\in H$, that $a(\act{\overline{2}}{h})=a(h^2)$, 
whence $(\overline{2}a)h=ah \cdot ah=(ah)^2$. 
Consequently, $\overline{2}h=(\overline{1}h)^2$ and $\overline{4}h=(\overline{2}h)^2=(\overline{1}h)^4$. 
Since $h$ is not the identity element of the group $N^H$, we deduce that $\overline{1}h\neq (\overline{0},\overline{1})$, 
the identity element of $N$. 
Since $N$ is a cyclic group of order $7$, we have $\overline{1}h\neq \overline{2}h$, $\overline{1}h\neq \overline{4}h$ and 
$\overline{2}h\neq \overline{4}h$. 
This means that $h$ is injective, and its image does not contain $(\overline{0},\overline{1})$.

{\it Case 2\/}: $(\overline{0},\overline{4})\iota_{11}=(t,\overline{2})$. 
A similar argument shows that $\overline{2}h=(\overline{1}h)^4$ and $\overline{4}h=(\overline{1}h)^2$, 
and we again deduce that $h$ is injective, and its image does not contain $(\overline{0},\overline{1})$.

By \eqref{D} and \eqref{E}, we have
\[ f_{p_{22}}^{11}=f_{(\overline{0},\overline{1})}^{12}p_{\xi_2\eta_2}^{H}f_{(\overline{0},\overline{1})}^{21} = 
(p_{\xi_2\eta_1}^{H})^{-1}p_{\xi_2\eta_2}^{H}(p_{\xi_1\eta_2}^{H})^{-1}, \]
and so
\begin{equation}\label{szendvicselem}
h=p_{\xi_1\eta_1}^H f_{p_{22}}^{11} = p_{\xi_1\eta_1}^H(p_{\xi_2\eta_1}^{H})^{-1}p_{\xi_2\eta_2}^{H}(p_{\xi_1\eta_2}^{H})^{-1}.
\end{equation}

This means that we can express $h$ as a product of entries in $P^H$ and their inverses, which sit at the intersections of two rows and 
two columns. 
By the definition of $P^H$, for any $a\in H$ and for any $i,j\in \{1,2\}$, we have
\[ ap_{\xi_i\eta_j}^H=p_{a\xi_i,a\eta_j}. \]
Hence the image of each entry of $P^H$ is contained in $\{(\overline{0},\overline{1}),p_{22}\}$, and
\[ ap_{\xi_i\eta_i}^H=p_{22}\textrm{ if and only if } a\xi_i=a\eta_j=2.\]
Consequently, for any $a\in H$,
\[ ap_{\xi_1\eta_1}^H=ap_{\xi_2\eta_2}^{H}=p_{22}\textrm{ if and only if } ap_{\xi_2\eta_1}^{H}=ap_{\xi_1\eta_2}^{H}=p_{22}. \]
Hence we see that it is impossible that two of the four entries 
$p_{\xi_1\eta_1}^H, p_{\xi_2\eta_1}^{H},$ $p_{\xi_2\eta_2}^{H}, p_{\xi_1\eta_2}^{H}$
sitting neither in the same row nor in the same column assign $p_{22}$ to some $a\in H$. 
For, in this case, \eqref{szendvicselem} would imply $ah=(\overline{0},\overline{1})$, contradicting
the property deduced above that the image of $h$ does not contain $(\overline{0},\overline{1})$.  
Consequently, for any $a\in H$, at least two of the four entries 
$p_{\xi_1\eta_1}^H, p_{\xi_2\eta_1}^{H}, p_{\xi_2\eta_2}^{H}, p_{\xi_1\eta_2}^{H}$
assign $(\overline{0},\overline{1})$ to $a$, and if precisely two then the respective entries sit
either in the same row or in the same column of $P^H$. 
So, by \eqref{szendvicselem}, we have $ah\in\{(\overline{0},\overline{1}),p_{22},p_{22}^{-1}\}$ for any $a\in H$, 
contradicting the fact that $ah\neq (\overline{0},\overline{1})$ and $h$ is injective.
This completes the proof that there is no embedding \eqref{beagyazas}
in the case of $S$ considered, thus proving the following result.

\begin{thm}
There exists a completely simple semigroup which is an extension of a completely simple semigroup $U$ by a group $H$ and 
which is not embeddable in the wreath product of $U$ by $H$.
\end{thm}

\section{Embeddability in a semidirect product}

In the previous section, we established that the Kaloujnine--Krasner Theorem does not generalize
for extensions of completely simple semigroups by groups. 
In this section, we present a modified version of the Kaloujnine--Krasner Theorem which holds for all extensions of 
completely simple semigroups by groups.

Let $S$ be an extension of a completely simple semigroup $U$ by a group $H$. 
Our goal is to give an embedding of $S$ into a semidirect product $V\rtimes H$ of a completely simple semigroup $V$ by $H$ 
such that, in the special case where $S$ is a group (i.e., $I$ and $\Lambda$ are singletons), it is just the embedding 
in \eqref{KKemb}. 
Unlike in the wreath product $U\wr H$, in this semidirect product $V\rtimes H$ the $\mathcal{R}$- and $\mathcal{L}$-classes of $V$, 
its sandwich matrix and the action of $H$ on $V$ can be chosen appropriately.

\begin{thm}
Any extension of a completely simple semigroup $U$ by a group $H$ is embeddable in a semidirect product of a 
completely simple semigroup $V$ by the group $H$, where the maximal subgroups of $V$ are direct powers of the maximal subgroups of $U$.
\end{thm}

\begin{biz}
Let $S$ be an extension of $U$ by $H$. 
As above, we can assume that $S=\mathcal{M}[G;I,\Lambda;P]$ where the sandwich matrix $P$ is normalized, 
and by Proposition \ref{kongruencia}, there is a normal subgroup $N$ of $G$ such that every entry of $P$ belongs to $N$, 
and $H=G/N$, $U=\mathcal{M}[N;I,\Lambda;P]\subseteq S$. 
Consider the action of $H$ on $N^H$ defining the wreath product $N\wr H$, and, 
for any $g\in G$, the map $f_g\in N^H$ defined in \eqref{KKemb}.

By means of $S$, we define a suitable semigroup $V$, an action of $H$ on $V$, 
and an embedding of $S$ into the semidirect product of $V$ by $H$. 
Let $V=\mathcal{M}[N^H;I,H\times \Lambda;Q]$, where the entries of $Q$ belong to the direct power $N^H$: 
for any $(B,\lambda)\in H\times \Lambda$ and $j\in I$, let
\[ q_{(B,\lambda),j}=\act{B}{f_{p_{\lambda j}}}. \]
Define an action of $H$ on $H\times \Lambda$ by the rule 
$\acth{A}{(B,\lambda)}=(A\cdot B,\lambda)$ $((B,\lambda)\in H\times \Lambda,\ A\in H)$. 
Now we give an action of $H$ on $V$ as follows: for any $A\in H$ and $(i,f,(B,\lambda))\in V$, let
\[  \acth{A}{(i,f,(B,\lambda))} = (i,\act{A}{f},\act{A}{(B,\lambda)}). \]
For any $A\in H$ and $(i,f,(B,\lambda)),(j,f',(C,\mu))\in V$, we have
\begin{eqnarray*}
\acth{A}{(i,f,(B,\lambda))}\cdot \acth{A}{(j,f',(C,\mu))} & = & (i,\act{A}{f},(A\cdot B,\lambda))(j,\act{A}{f'},(A\cdot C,\mu)) \\
& = & (i,\act{A}{f}\cdot q_{(A\cdot B,\lambda),j}\cdot \act{A}{f'},(A\cdot C,\mu))\\
& = & (i,\act{A}{f}\cdot \act{A}{q_{(B,\lambda),j}}\cdot \act{A}{f'},\acth{A}{(C,\mu)})\\
& = & \acth{A}{(i,f\cdot q_{(B,\lambda),j}\cdot f',(C,\mu))}\\
& = & \acth{A}{\big( (i,f,(B,\lambda))(j,f',(C,\mu)) \big)}.
\end{eqnarray*}
Hence this is a well-defined action of $H$ on $V$, and so the semidirect product 
$V\rtimes H=\mathcal{M}[N^H;I,H\times \Lambda;Q]\rtimes H$ 
with respect to this action is defined.

Let us consider the mapping
\[ \psi\colon \mathcal{M}[G;I,\Lambda;P] \rightarrow \mathcal{M}[N^H;I,H\times \Lambda;Q]\rtimes H, \]
where
\[ (i,g,\lambda)\psi = ((i,f_g,(gN,\lambda)),gN). \]

We intend to verify that $\psi$ is an embedding. 
Assume that $(i,g,\lambda)\psi=(j,h,\mu)\psi$, i.e., $(i,f_g,(gN,\lambda)),gN)=(j,f_h,(hN,\mu)),hN)$. 
Hence $i=j$, $\lambda=\mu$, $gN=hN$ and $f_g=f_h$. 
Since $\kappa$ in \eqref{KKemb} is injective, the last two equalities imply $g=h$, and so $(i,g,\lambda)=(j,h,\mu)$ follows.

To prove that $\psi$ is a morphism, we can see for any $(i,g,\lambda),(j,h,\mu)\in \mathcal{M}[G;I,\Lambda;P]$, that
\begin{eqnarray*}
(i,g,\lambda)\psi(j,h,\mu)\psi & = & ((i,f_g,(gN,\lambda)),gN)((j,f_h,(hN,\mu)),hN)\\
& = & ((i,f_g,(gN,\lambda))\cdot \act{gN}{(j,f_h,(hN,\mu))},ghN)\\
& = & ((i,f_g,(gN,\lambda))(j,\act{gN}{f_h},(ghN,\mu)),ghN)\\
& = & ((i,f_g\cdot q_{(gN,\lambda),j}\cdot \act{gN}{f_h},(ghN,\mu)),ghN),
\end{eqnarray*}
and
\[ \big( (i,g,\lambda)(j,h,\mu)\big)\psi=(i,gp_{\lambda j}h,\mu)\psi =  ((i,f_{gp_{\lambda j}h},(gp_{\lambda j}hN,\mu)),gp_{\lambda j}hN).\]
We need to prove that the two maps in the middle components are equal. 
Since $p_{\lambda j} \in N$ and $N$ is the identity element of $H$, \eqref{hatas} implies by the definition of $Q$ that
\begin{eqnarray*}
f_{gp_{\lambda j}h} & = & f_g\cdot \act{gN}{f_{p_{\lambda j}h}}\\
& = & f_g\cdot \act{gN}{(f_{p_{\lambda j}}\cdot \act{p_{\lambda j}N}{f_h})}\\
& = & f_g\cdot \act{gN}{(f_{p_{\lambda j}}\cdot f_h)}\\
& = & f_g\cdot \act{gN}{f_{p_{\lambda j}}} \act{gN}{f_h}\\
& = & f_g\cdot q_{(gN,\lambda),j}\cdot \act{gN}{f_h}.
\end{eqnarray*}
Thus $\psi$ is, indeed, an embedding, and the proof of the theorem is complete.

\end{biz}

Note that, in the case where $S$ is a group, i.e., where $I$ and $\Lambda$ are singletons 
(and so the single entry of $P$ is the identity of $G$, and $S$ is isomorphic to $G$), 
the map $\psi$ coincides with the embedding $\kappa$ in \eqref{KKemb}.

\end{document}